\RequirePackage[english]{babel}
\documentclass{amsart}

\usepackage{amssymb,amscd,amsthm}

\DeclareMathOperator{\Tt}{\mathbb T}
\DeclareMathOperator{\R}{\mathbb R}

\DeclareMathOperator{\Z}{\mathbb Z}

\DeclareMathOperator{\dom}{dom}

\def\rf{{{\mathbb R}_{\text{\scriptsize field}}}}

\theoremstyle{plain}
\newtheorem{theorem}{Theorem}
\newtheorem{lemma}[theorem]{Lemma}

\newtheorem{proposition}[theorem]{Proposition}
\newtheorem{corollary}[theorem]{Corollary}

\newtheorem{fact}[theorem]{Fact}
\theoremstyle{definition}
\newtheorem{remark}[theorem]{Remark}
\newtheorem{definition}[theorem]{Definition}
\newtheorem{example}[theorem]{Example}
\newtheorem{exercise}[theorem]{Exercise}

\numberwithin{theorem}{section}

\newcommand{\bt}{\begin{theorem}}
\newcommand{\et}{\end{theorem}}

\newcommand{\bl}{\begin{lemma}}
\newcommand{\el}{\end{lemma}}

\newcommand{\bfa}{\begin{fact}}
\newcommand{\efa}{\end{fact}}

\newcommand{\bexa}{\begin{example}}
\newcommand{\eexa}{\end{example}}

\newcommand{\bexe}{\begin{exercise}}
\newcommand{\eexe}{\end{exercise}}

\newcommand{\bprop}{\begin{proposition}}
\newcommand{\eprop}{\end{proposition}}

\newcommand{\bp}{\begin{proof}}
\newcommand{\ep}{\end{proof}}

\newcommand{\bc}{\begin{corollary}}
\newcommand{\ec}{\end{corollary}}

\newcommand{\bd}{\begin{definition}}
\newcommand{\ed}{\end{definition}}

\newcommand{\br}{\begin{remark}}
\newcommand{\er}{\end{remark}}

\def\MM#1{\raise0.7ex\hbox{\tiny\ding{170}}\marginpar{$\mathcal{M}^2$:
{\footnotesize #1}}}

\title{Topology of definable abelian groups in o-minimal structures}

\author{El\'ias Baro}
\address{Departamento de \'Algebra,
Facultad de Matem\'aticas,
Universidad Complutense de Madrid,
28040 Madrid}
\email{\tt eliasbaro@pdi.ucm.es}

\author{Alessandro Berarducci} 
\address{Universit\`a di Pisa, Dipartimento di Matematica, Largo Bruno
Pontecorvo 5, 56127 Pisa, Italy}
\email{\tt berardu@dm.unipi.it}
\thanks{The first author is partially supported
by MTM2008-00272 and
Grupos UCM 910444. The second author is partially supported by PRIN
2007PRYAAF 004: O-minimalit\`a - Metodi e modelli non
standard - Teoria degli insiemi.}



\date{Feb. 11, 2011}

\begin{document} 
\begin{abstract}In this note we prove that every definably connected,
definably compact abelian definable group $G$ in an o-minimal expansion of a
real closed field with $\dim(G)\neq 4$ is definably homeomorphic to a torus of
the same dimension. Moreover, in the semialgebraic case the result holds for
all dimensions.
\end{abstract}

\maketitle


\section{Introduction} Let $M$ be an o-minimal expansion of a real closed field.
Let $H$ be a definable group in $M$ equipped with Pillay's topology in
\cite{Pi:88}. So when
$M$ is an expansion of the real line $H$ is a real Lie group, and in general it
is an ``$M$-Lie group''. An example of definable group is the {\em n-torus}
$\Tt^n(M)$ over $M$, defined as the poly-interval $[0,1)^n$ in $M$ with the sum
operation modulo $1$. When $M$ is the real field $\rf$ this coincides with the
classical torus $(\R/\Z)^n$. Let us recall that $(\R/\Z)^n$ is the only compact
abelian connected Lie group of dimension $n$ up to Lie-isomorphisms.

In the rest of the paper we fix a definably connected, definably
compact, definable abelian group $G$ in $M$ of dimension $n$, endowed with
Pillay's topology. There are three natural questions:

\begin{enumerate}
\item  Is $G$ definably isomorphic to the $n$-torus $\Tt^n(M)$?
\item Is $G$ definably homeomorphic to $\Tt^n(M)$?
\item  Is $G$ definably homotopy equivalent to $\Tt^n(M)$?
\end{enumerate}

The answer to question (1) is clearly no. For instance $[0,1)\subset \mathbb{R}$
 modulo $1$ is Lie isomorphic to $SO(2,\R)$ but the isomorphism is not
semialgebraic (note, however, that they are semialgebraically
homeomorphic). Instead, we could ask if $G$ is definably isomorphic to a
product of $1$-dimensional definable subgroups. But it turns out that this
question still has a negative answer even for $M = \rf$. Indeed it is
possible that $\dim(G)>1$, but $G$ has no subgroups of dimension one definable
in $\rf$
\cite[Example 5.2]{PeSte:99}.

On the other hand the homotopy problem (3) has a positive answer by
\cite[Theorem 3.4]{BeMaOt:08}. 

In this note we deal with the homeomorphism problem (2), giving a positive
solution when $\dim(G)\neq 4$ (see Theorem \ref{deftheorem}). We also show that
when $G$ is semialgebraic, namely $M$ is a real closed field without additional
structure, then the
result holds in all dimensions (see Theorem \ref{semialgtheorem}). We point out
that the one dimensional case was already proved in \cite{St:94}.

The assumption that $G$ is definably compact and definably connected is not very
restrictive.
In fact by \cite[Thm.5.1, Thm.5.7]{PeSta:05} and
\cite[Thm.1.2]{PeSte:99}, the study of the topology of definable abelian
groups can be reduced to the definably compact case. 
So as a corollary we obtain, with a small proviso, that a definably connected abelian
n-dimensional group is definably homeomorphic to a space of the form
$\Tt^m(M)\times M^k$ with $m+k = n$. The proviso is vacuous in the semialgebraic
case, while in the general case we require that $m\neq 4$, where $m$ is $n$
minus the dimension of the maximal torsion free definable subgroup of $G$.

Suitable versions of problems (1),(2),(3) can be posed in the non-abelian case.
This can be done as follows. Given a definable group $G$, there is a canonical
real Lie group $G/G^{00}$ associated to $G$ (by \cite{Pi:04} and
\cite{BOPP:05}). By \cite{Ba:10} and \cite{BeMa:10} when $G$ is definably
compact and definably connected, the isomorphism type of $G/G^{00}$ determines
$G$ up to definable homotopy equivalence. One can ask whether the isomorphism
type of $G/G^{00}$ determines $G$ up to definable homeomorphism. The results in
\cite{Ma:10} reduce the question to the abelian case, which is the one
we consider in this paper. Finally let us observe that by
\cite[Thm.3.8.8]{Co:09t} (see also \cite{Co:09}) the study of the topology of a
definable group reduces to the
definably compact case.

We shall make use of the ``o-minimal Hauptvermutung'' proved by M. Shiota in
\cite[Chapter III]{Sh:97} when $M$ is an expansion of $\rf$, and extended in
\cite[\S 2]{Sh:10} to the case when $M$ is an o-minimal expansion of
an arbitrary real closed field. 

\bfa\label{haupt}\emph{(o-minimal Hauptvermutung)} Let $K$ and $L$ be
finite simplicial complexes. Let $M$ be an o-minimal expansion of a real closed field.  
If there exists a definable, in $M$, homeomorphism from
$|K|$ to $|L|$, then there is a $PL$-isomorphism from $|K|$ to $|L|$.
\efa

Here $|K|$ denotes a geometrical realization, in $M$, of the simplicial complex
$K$.
In this note all simplicial complexes are closed and finite, so that $|K|(M)$ is always definably compact. By a ``PL map'' we always mean ``finitely PL map'', namely the geometrical realization of a simplicial map between finite subdivisions of the relevant complexes. In Section 2 we prove the semialgebraic case of the homeomorphism problem (2) using the Hauptvermutung in the weak form of \cite{Sh:97}. In Section 3 we assume $\dim(G)\neq 4$ and we reduce the general o-minimal case to the semialgebraic case. In this step we need the strong form of the Hauptvermutung (as in \cite{Sh:10}) and the following fact: 

\bfa \label{classtori}\emph{(Classification of Homotopy tori)} Let $X$ be a
closed
$PL$-manifold
of dimension
$n\neq 4$ homotopy equivalent to the standard torus $\mathbb{T}^n(\R)$ (considered as a $PL$-manifold under a standard triangulation). Then
there is a finite covering $f:\widetilde{X}\rightarrow X$ such that
$\widetilde{X}$ is $PL$-homeomorphic to $\mathbb{T}^n(\R)$.
\efa

When $n\geq 5$ a proof of Fact \ref{classtori} can be found in \cite[Theorem
B]{HsSi:69} and \cite[Corollary]{Wa:69}. (See also \cite[Chapter 15A]{Wa:99} for
a complete development of homotopy tori.) When $n\leq 3$ it turns out that $X$
is already $PL$-homeomorphic to a standard torus. Indeed, for $\dim(X)=1$ or $2$
this is well-known and for $\dim(X)=3$ we can use \cite[Theorem 5.4, pag
249]{KiSi:77} together with the positive solution of the three dimensional
Poincar\'e's conjecture. Since our intended readership may not be familiar with
the notations in \cite{KiSi:77} we add few lines of explanation. The cited
theorem tells us that $S^*(\Tt^3) = 0$. Unraveling the notations this means that
if $M_1$ and $M_2$ are $PL$-manifolds and $f_i\colon M_i \to \Tt^3$ is a
homotopy equivalence for $i=1,2$, then there is a $PL$-homeomorphism $h\colon
M_1\to M_2$ such that $f_2 \circ h$ is homotopic to $f_1$.
In particular, taking $M_2 = \Tt^3$ and $f_2 = id$, we obtain that any
$PL$-manifold homotopy equivalent to $\Tt^3$ is $PL$-homeomorphic to $\Tt^3$.
This statement (= Borel's conjecture for $\Tt^3$) is known to
imply
Poincar\'e's conjecture in dimension 3 \cite[\S 1.4]{F:96}, which was not known
when \cite{KiSi:77} was written. The solution of the riddle lies in a note
hidden inside the proof of Theorem 5.3 in \cite{KiSi:77} whose effect is to
modify the definition of $S^*$ in dimension 3: ``in dimension $3$ we supplement
this definition by supposing that $M_1$ is Poincar\'e, i.e., contains no fake
3-discs''. Granted the positive solution to the 3-dimensional Poincar\'e's
conjecture, the supplement is vacuous.

\section{Semialgebraic case}

In this section suppose that $M$ is a real closed field without additional structure. So the definable sets in $M$ coincide with the semialgebraic sets. We prove. 

\bt \label{semialgtheorem}
$G$ is semialgebraically homeomorphic to the
$n$-standard torus $\Tt^n(M)$.
\et

\bp By Robson's embedding theorem (see \cite[Theorem 10.1.8]{vdD:98}) we can
assume that the topology of $G$ (given by \cite{Pi:88}) coincides with the
topology induced by the ambient space $M^n$. By the triangulation theorem we can
then assume that the underlying set $\dom(G)$ is the realization of a
$\emptyset$-definable finite simplicial complex $K$. A priori we cannot ensure
that the group operation of $G$ is $\emptyset$-definable, but  by model
completeness of the theory of real closed fields there exist a possibly
different group operation on $\dom(G)=|K|$ which is $\emptyset$-definable and
continuous with respect to the topology of $|K|$. Since
we are only interested in the definable homeomorphism type of $G$ we can assume
the group operation is $\emptyset$-definable. We can then consider the group
$G(\R)$ obtained by interpreting the defining formulas in $\rf$. By
\cite[Remark 2.6]{Pi:88}, there is a (unique) Nash group structure on $G(\R)$.
In particular, $G(\R)$ is an abelian compact connected real Lie-group and
therefore there is a Lie-isomorphism $f\colon G(\R)\to \Tt^n(\R)$. We will show
that $f$ is definable in some o-minimal expansion of the real field. In fact it
is enough to consider the o-minimal structure $\R_{an}$ studied in
\cite{vdD:86}. We need the following:

{\em Fact}: Given an analytic function $f$ defined on an open subset $V$ of $\R^n$, its restriction to a definable (i.e. semialgebraic) compact subset $K\subset V$ is definable in $\R_{an}$. 

Indeed this is true (almost by definition of $\R_{an}$) when $K$ is a compact poly-interval, and the general case follows by covering $K$ by finitely many poly-intervals contained in $V$. We then obtain:

{\em Claim}: The Lie-isomorphism $f\colon G(\R)\to \Tt^n(\R)$ is definable in $\R_{an}$. 

In fact there are semialgebraic charts making $G(\R)$ into a Nash group and for
each chart $V$, $f|V$ is analytic. By shrinking the charts we can assume that
$f|V$ extends to an analytic map on the closure of $V$. So by the above fact $f$
is definable in $\R_{an}$.

In particular we have proved that there is a homeomorphism $f\colon G(\R)\to \Tt^n(\R)$ definable in $\R_{an}$. By the semialgebraic triangulation theorem and the o-minimal Hauptvermutung of \cite{Sh:97}, there is a semialgebraic homeomorphism $g\colon G(\R)\to \Tt^n(\R)$. Moreover, 
by model completeness of the theory of real closed field, there is some $g$ as above which is $\emptyset$-definable in $\rf$. Interpreting the same formulas in $M$ we obtain a semialgebraic homeomorphism from $G(M)$ to $\Tt^n(M)$ as desired. \ep

\section{General o-minimal case}
In this section we assume that $M$ is an arbitrary o-minimal expansion of a real closed field. 
We will prove: 

\bt \label{deftheorem} If $n=\dim(G) \neq 4$, $G$ is definably
homeomorphic to the $n$-torus $\mathbb{T}^n(M)$. 
\et

As above, we can assume that Pillay's topology on $G$ coincides with the 
topology induced by the ambient space $M^n$, and by the triangulation theorem we can then assume that $\dom(G)$ is the geometrical realization $|K|(M)$ of a finite simplicial complex $K$. We need: 

\bl \label{lemma} If $n=\dim(G) \neq 4$, $\dom(G)=|K|(M)$ admits a semialgebraic
abelian group operation (possibly unrelated to the original one). \el

Theorem \ref{deftheorem} follows at once from the Lemma and the semialgebraic
case (Theorem \ref{semialgtheorem}). So it remains to prove the lemma.

\bp[Proof of Lemma \ref{lemma}] Note that $\dom(G) = |K|(M)$ is at the same time
a closed definable manifold (with Pillay's topology) and the realization, over
$M$, of a finite simplicial complex. By
Shiota's o-minimal Hauptvermutung in \cite{Sh:10}, it easily follows (see Fact \ref{PLmanifold} below) that
$|K|(M)$ is a closed $PL$-manifold ``over $M$''. This is equivalent to say that
the closed star of each
vertex of $K$ is $PL$-homeomorphic to the standard simplex of the correct
dimension. By model completeness of the theory of real closed fields, the same
holds over $\R$. Namely $|K|(\R)$ is a $PL$-manifold (but we have no way of
inheriting the definable group structure of $|K|(M)$). Moreover $|K|(\R)$ is
homotopy equivalent to the standard torus
$\Tt^n(\R)$. Indeed, by \cite[Theorem 3.4]{BeMaOt:08} there exists a definable
homotopy equivalence from $\dom(G) = |K|(M)$ to $\Tt^n(M)$ and therefore by \cite[Theorem
3.1]{BaOt:10} there is a semialgebraic homotopy equivalence from $|K|(\R)$ to
$\Tt^n(\R)$.
Since $n=\dim(G)\neq 4$, by Fact \ref{classtori},
$|K|(\R)$ has a finite $PL$-cover which is $PL$-homeomorphic to $\Tt^n(\R)$. 
Namely we have a PL-covering $f:\Tt^n(\R)\to |K|(\R)$ with finite fibers. 
By model completeness we can assume that $f$ is defined without parameters. So
we can go back to $M$ and get a semialgebraic (actually $PL$) covering $$f\colon
\Tt(M)^n\to |K|(M) = \dom(G).$$ 
But on $\dom(G)$ we have a
definable group operation that can be lifted to
$\Tt^n(M)$ via $f$ (by uniform lifting of paths). 
So we get a definable group operation $*$ on $\Tt^n(M)$ making $f$ into a definable covering homomorphism with a finite kernel $\Gamma < (\Tt^n(M),*)$. Note that $*$ may not coincide with the natural group operation on $\Tt^n(M)$ (the sum mod 1), so in particular it need not be semialgebraic. 
In any case however $(\Tt^n(M),*)$ is an abelian group. Therefore there
is $k$ such that $\Gamma$ is contained in the $k$-torsion subgroup
$(\Tt^n(M),*)[k]$ of $(\Tt^n(M),*)$.
Our next goal is to obtain a definable group homomorphism $$h \colon G \to (\Tt(M)^n,*).$$
Write for simplicity $\Tt^n$ for $\Tt^n(M)$. Now $G$ is definably isomorphic to $(\Tt^n,*)/\Gamma$ and since $\Gamma < (\Tt^n,*)[k]$ 
there is a definable covering from $(\Tt^n,*)/\Gamma$ to $(\Tt^n,*)/(\Tt^n,*)[k]$. The latter group 
is definably isomorphic to $(\Tt^n,*)$ because the $k$-torsion subgroup of a
definably compact definably connected abelian group $H$ is the kernel of the
surjective homomorphism $H \to H$ sending $x$ to $kx$ (we need the fact that
such groups are divisible, as proved in \cite[Theorem 2.1]{EdOt:04}). 
Composing we obtain a finite definable covering homomorphism $h\colon G\to (\Tt^n,*)$. 
As already remarked $*$ may not be semialgebraic. However on $\Tt^n$ we also
have a semialgebraic group operation $\cdot$ (the sum mod 1).  The idea is to
use the covering map $h$ (seen just as a continuous map, not as a group
homomorphism) to lift the semialgebraic group operation $\cdot$ to a
semialgebraic group operation on $\dom(G)$. The problem however is that $h$ is
not semialgebraic. However by \cite[Corollary 2.2]{EdJoPe:10}, each definable
cover of a semialgebraic group, is equivalent to a semialgebraic cover. So there
is a semialgebraic covering homomorphism
$h'\colon G'\to (\Tt^n,\cdot)$ and a definable homeomorphism
$\psi:\dom(G)\to \dom(G')$ commuting with $h$ and $h'$. But $\dom(G)$ and $\dom(G')$ are
semialgebraic, so by the Hauptvermutung (combined with the triangulation theorem) 
there is a semialgebraic homeomorphism $\phi:\dom(G)\to \dom(G')$. Now take $h'\circ \phi$. This is a semialgebraic covering from $\dom(G)$ to $(\Tt^n,\cdot)$, and it can be used to lift $\cdot$ to a semialgebraic group operation on $\dom(G)$. 
\ep 

Let us prove the missing fact needed in the above proof.

\bfa \label{PLmanifold} Let $K$ be a finite simplicial complex
such that $|K|$ is an
$n$-dimensional closed definable manifold. Then $|K|$ is a PL-manifold, namely the star
of each vertex of $K$ is $PL$-homeomorphic to the standard $n$-simplex.
\efa
\bp In this proof simplicial complexes are assumed to
be finite but not necessarily closed. We will use without mention the well-known
invariance of stars in
piecewise linear topology, i.e., the star of a vertex of a closed simplicial
complex is
PL-isomorphic to the star of that vertex in any simplicial subdivision.
Let $\{(U_1,f_1),\ldots,(U_s,f_s)\}$ be a definable atlas of
$|K|$. That is, each $U_i$ is a definable open subset of
$|K|$, each $f_i$ is a definable homeomorphism from
$U_i$ to a definable open subset $V_i$ of $M^n$ (with the usual property
on transition maps) and $|K|=\bigcup_{i=1}^s U_i$. By shrinking of
coverings, we can find
definable open subsets $W_1,\ldots,W_s$ of $|K|$ such that
$|K|=\bigcup_{i=1}^s W_i$ and $W_i\subset
\overline{W}_i\subset U_i$ for each $i=1,\ldots,s$. Moreover, by the
triangulation theorem we
can assume that each $V_i$ is the realization of an open finite simplicial
complex and
$f_i(\overline{W_i})$ the realization of a subcomplex. Considering a
barycentric subdivision if necessary, we can also assume that the star in
$V_i$ of each
vertex of $f_i(\overline{W_i})$ is a closed finite subcomplex. In
particular, since $V_i$ is an open subset of $M^n$, it follows that the
star of each vertex in
$f_i(\overline{W_i})$ is PL-isomorphic to a standard $n$-simplex.

Now, again by the triangulation theorem, there exist a definable homeomorphism
$\psi:|L|\rightarrow |K|$ compatible with the definable sets $U_i$, $W_i$ and
$\overline{W_i}$.  Since $\psi^{-1}(\overline{W_i})$ and
$f_i(\overline{W_i})$ are definable homeomorphic, by the o-minimal
Hauptvermutung they are PL-isomorphic. Now, given a vertex $v$ of $L$, the star
of $v$ in $L$ is contained in some $\overline{W_i}$ and therefore is
PL-isomorphic to the star of some vertex of $f_i(\overline{W_i})$. Hence, the
star of each vertex of $L$ is PL-isomorphic to a standard $n$-simplex.

By the o-minimal Hauptvermutung, there exist a PL-isomorphism of $|L|$
and $|K|$. Hence
we deduce that the star of each vertex of $K$ is PL-isomorphic to the star of a
vertex of $L$, which in turn is PL-isomorphic to a standard $n$-simplex.
\ep

We have thus completed the proof of the Lemma, and Theorem \ref{deftheorem} follows.  

A possible attempt to deal with the case $\dim(G)=4$ is to replace $G$ with
$G\times \Tt^1$. So by Theorem \ref{deftheorem} there is
definable homeomorphism from $G \times \Tt^1$ to $\Tt^{5}$. However we do not
know whether this implies that there is a definable homeomorphism from $G$ to
$\Tt^4$. Finally let us observe that, even in dimension $4$,  we can always
assume $\dom(G)=|K|(M)$ (after a triangulation) and conclude that $|K|(\R)$ is
homotopy equivalent to $\Tt^4(\R)$ (reasoning as in the proof of Lemma
\ref{lemma}). Moreover by \cite[\S 11.5]{FrQu:90} a fourth dimensional
$PL$-manifold homotopy
equivalent to
a standard torus is homeomorphic to it (although not necessarily
$PL$-homeomorphic). So in any case we conclude that $|K|(\R)$ is homeomorphic to
$\Tt^4(\R)$, but since a priori the homeomorphism could be quite wild, there is
no obvious way to obtain from these data a definable homeomorphism from $|K|(M)$
to $\Tt^4(M)$.

\section*{Acknowledgements} We thank M. Mamino for many discussions regarding the topics of this paper, 
and R. Benedetti, R. Frigerio, A.
Ranicki, M. Shiota and C.T.C. Wall for their comments concerning low
dimensional topology and classification of manifolds. This work was done while the first author was visiting the
Department of Mathematics of the University of Pisa in the period 6 January - 27 July 2010. 

\begin{footnotesize}
\end{footnotesize}
\end{document}